\documentclass{ifacconf}
\usepackage{amsmath}
\usepackage{amsfonts}
\usepackage{amssymb}
\usepackage{mathtools}
\usepackage{tikz}
\usepackage{xcolor, color, soul}
    \sethlcolor{green}
    \definecolor{hbl}{rgb}{0.0, 0.48, 0.65}
\usepackage{graphicx}      
\usepackage{natbib}

\newcommand{\R}{\mathbb{R}}
\newcommand{\N}{\mathbb{N}}

\newcommand{\tr}{\operatorname{tr}}

\newcommand{\Real}{\operatorname{Re}}

\newcommand{\drh}[2]{\Omega^{#2}(#1)}

\DeclarePairedDelimiterX{\norm}[1]{\lVert}{\rVert}{#1}
\DeclarePairedDelimiterX{\abs}[1]{\lvert}{\rvert}{#1}
\DeclarePairedDelimiterX{\scprod}[2]{\langle}{\rangle}{#1, #2}
\DeclarePairedDelimiterX{\dualprod}[2]{\langle}{\rangle}{#1, #2}

\begin{document}
\begin{frontmatter}

\title{Well-Posedness of the Hodge Wave Equation on a Compact Manifold \thanksref{footnoteinfo}} 

\thanks[footnoteinfo]{Funded by the European Union (Horizon Europe MSCA project ModConFlex, grant number
101073558). Views and opinions expressed are however those of the author only and do
not necessarily reflect those of the European Union or the European Research Council
Executive Agency. Neither the European Union nor the granting authority can be held
responsible for them.}

\author[First]{Filippo Testa} 

\address[First]{Department of Applied Mathematics, University of Twente, 
   Enschede, The Netherlands (e-mail: f.testa@utwente.nl).}

\begin{abstract}           
    In this work, we study the Hodge wave equation on a compact orientable manifold. We present the necessary differential geometry language to treat Sobolev spaces of differential forms and use these tools to identify a boundary triplet for the problem.
    We use this boundary triplet to determine a class of boundary conditions for which the problem is well-posed.
\end{abstract}

\begin{keyword}
Wave equation, manifold, external calculus, differential forms, port-Hamiltonian, boundary triplet
\end{keyword}

\end{frontmatter}
\section{Introduction}
The ubiquitous presence of the wave equation in applications, as well as its fundamental role as the main example and gateway for the broader class of hyperbolic equations, provide the main motivation behind this paper.\\
As a simple yet powerful tool for modeling small oscillations, it is natural to formulate the wave equation on a manifold representing the oscillating object. This manifold is often not guaranteed to simply be a bounded domain in the Euclidean space; for example, when modeling membranes, shells, or airplane wings, we may be interested in using a curved, two-dimensional manifold.\\
The language of differential geometry is then necessary to treat the problem, as it takes care of curvature in an elegant and effective way. This necessity, which may turn us away at first, rewards with a global approach to the problem, which better incorporates the geometry of the object into the study of its oscillations, revealing insights that would be lost otherwise.\\
In this paper, we investigate the well-posedness of the wave equation on a compact, orientable, Riemannian manifold $M$, which we refer to as the Hodge wave equation.
In Section $2$, we introduce the necessary differential geometry language, first recalling the basics for smooth differential forms, then expanding the language to $L^2$ forms and Sobolev spaces. In particular, we derive an extension of the Stokes theorem for $L^2$ differential forms, which, in Section $4$ allows us to turn the question of well-posedness of the Hodge wave equation on a manifold into a question about boundary values, following the functional analytic theory presented in Section $3$.\\
We remark that the method we use for classifying boundary conditions is the standard one for port-Hamiltonian systems (see \cite{ScM:02}, \cite{JZ:12}) and was also used in \cite{KZ:15} for studying the wave equation on a domain in $\R^n$.

\section{Sobolev theory on compact orientable manifolds}
\subsection{Smooth differential forms}
In this section we recall some concepts from smooth differential geometry and focus on extending them to the $L^2$ setting. For more details on the smooth theory see \cite{Lee:12}, while for the Sobolev theory of differential forms we recommend \cite{AR:06}.\\
Consider a compact, orientable, $n$-dimensional smooth Riemannian manifold $M$ with boundary.\\
For every $k\in\N$ let $\drh{M}{k}$ be the space of smooth differential forms of degree $k$ on $M$. Upon picking a sufficiently small coordinate patch $U\subset M$,  a $k$-form $\omega$ may be transported onto an open subset $V$ of the closed half space $\R^n_+$. Assuming $x_i$ for $i=1,\dots,n$ are the coordinates on $V$, we can then write $\omega$ (using Einstein's summation convention) as
\begin{equation}\label{coord}
    \omega = \omega_{i_1\dots i_k} dx^{i_1}\wedge \dots\wedge dx^{i_k}
\end{equation}
for some $\omega_{i_1\dots i_k} \in C^{\infty}$ which we will some times refer to as the coefficients of $\omega$ in this coordinate patch.\\
Spaces of forms of different degree are interconnected by the wedge product
\[
    \omega\in\drh{M}{k},\mu\in\drh{M}{h}\Rightarrow\omega\wedge\mu\in\drh{M}{k+h},
\]
which, on a coordinate patch, can be written as
\begin{align}\label{wedge}
    &(\omega_{i_1\dots i_k} dx^{i_1}\wedge \dots\wedge dx^{i_k})\wedge (\mu_{j_1\dots j_k} dx^{j_1}\wedge \dots\wedge dx^{j_k})=\nonumber\\
    & \omega_{i_1\dots i_k}\mu_{j_1\dots j_k}dx^{i_1}\wedge \dots\wedge dx^{i_k}\wedge dx^{j_1}\wedge \dots\wedge dx^{j_k}.
\end{align}
This easily results in the formula
\begin{equation}\label{swap}
    \omega\wedge\mu=(-1)^{hk}\mu\wedge\omega,\quad \forall\omega\in\drh{M}{k},\forall\mu\in\drh{M}{h}.
\end{equation}
For any $k\in \N$ we may define the external differential map
\[
    d\colon \drh{M}{k}\to\drh{M}{k+1}
\]
which, for a form $\omega$ on a coordinate patch as (\ref{coord}), has the explicit formula
\[
    d\omega = \frac{\partial\omega_{i_1\dots i_k}}{\partial x^j} dx^j\wedge dx^{i_1}\wedge \dots\wedge dx^{i_k}.
\]
This description allows us to obtain the following chain rule for $\omega\in\drh{M}{k},\mu\in\drh{M}{h}$:
\begin{equation}\label{chain}
    d(\omega\wedge\mu)=d\omega\wedge\mu+(-1)^k\omega\wedge d\mu.
\end{equation}
Assuming we have defined a positive scalar product on the tangent space of $V$, which, for two vector fields $X,X\in \mathfrak{X}(V)$ is given by
\[
    \scprod{X^i\frac{\partial}{\partial x^i}}{Y^j\frac{\partial}{\partial x^j}}_g\coloneqq X^ig_{ij}Y^j,
\]
where $g_{ij}$ are the parameters of the inner product for the basis $\{\frac{\partial}{\partial x^1},\dots,\frac{\partial}{\partial x^n}\}$. This scalar product then induces a scalar product onto the dual of the tangent space of $V$:
\[
    \scprod{\omega_i dx^i}{\mu_jdx^j}_g\coloneqq \omega_ig^{ij}\mu_j,
\]
where the matrix $[g^{ij}]_i^j$ is the inverse of $[g_{ij}]_i^j$, that is, $g^{ij}g_{jk}=g_{kj}g^{ji}=\delta_k^i$.\\
Finally, this scalar product induces a product onto the space of $k$-forms on $V$ (which are well-defined on open domains in $\R^n_+$):
\begin{align*}
    &\scprod{dx^{i_1}\wedge \dots\wedge dx^{i_k}}{dx^{j_1}\wedge \dots\wedge dx^{j_k}}_g\\
    &\coloneqq \sum_{\sigma\in S_k}(-1)^{\sigma}g^{i_1j_{\sigma(1)}}\dots g^{i_kj_{\sigma(k)}},
\end{align*}
where $S_k$ is the group of permutations of $k$ elements, and $(-1)^{\sigma}$ is the sign of the permutation $\sigma\in S_k$. More precisely, for $\omega^1,\omega^2\in \drh{V}{k}$, we have
\[
    \scprod{\omega^1}{\omega^2}=\omega_{i_1\dots i_k}^1\omega_{j_1\dots j_k}^2g^{i_1j_1}\dots g^{i_kj_k}.
\]
Since the manifold $M$ we are considering is Riemannian, that is, we have made a choice of the Riemannian metric $g$, there is a coherent choice of a positive scalar product on the tangent space at every point, hence the scalar product on differential forms is well defined on $\drh{M}{k}$.

\subsection{$L^2$ differential forms and Sobolev spaces}
The choice of a Riemannian metric guarantees the existence of a volume form $\operatorname{vol}\in \drh{M}{n}$ which may be integrated over the manifold. This allows us to define the space of square-integrable functions on the manifold $M$ as the closure of $C^{\infty}(M)$ with respect to the norm
\[
\norm{f}_{L^2}^2\coloneqq \int_M f^2\operatorname{vol}.
\]
We refer to this space as $L^2(M)$ or $L^2\drh{M}{0}$.\\
In order to extend the concept of $L^2$ to forms of any degree, we define a new inner product on $\drh{M}{k}$ for $k=0,\dots,n$ as
\begin{equation}\label{l2p}
    \scprod{\omega}{\mu}_{L^2\drh{M}{k}}\coloneqq \int_M \scprod{\omega}{\mu}_g\operatorname{vol}.
\end{equation}
We can define the Hodge star map
\begin{equation}
    \star\colon\drh{M}{k}\to\drh{M}{n-k},
\end{equation}
where $\star\mu$ is the (unique) $n-k$-form such that
\begin{equation}\label{HStar}
    \scprod{\omega}{\mu}_g\operatorname{vol}=\omega\wedge\star\mu \quad \forall\omega\in\drh{M}{k}.
\end{equation}
Notice that the term on the right is well-defined since $\omega$ is a $k$-form, while $\star\mu$ is a $n-k$-form, hence their wedge product is a $n$-form and can be integrated.\\
The Hodge star is an isometry with respect to the $L^2$ product, and the following useful formula holds:
\begin{equation}\label{invHS}
    \star^{-1}=(-1)^{(n-k)k}\star.
\end{equation}
Through the Hodge star we may define the formal adjoint of the external differential, the codifferential operator
\begin{equation}\label{delta}
    \delta\colon \drh{M}{k}\to\drh{M}{k-1},\quad \delta\coloneqq (-1)^k\star^{-1}d\star.
\end{equation}
The most important result about the external differential is the Stokes theorem, which states that for $\omega\in \drh{M}{n-1}$ it holds that
\begin{equation}\label{Stokes}
    \int_M d\omega = \int_{\partial M}\omega.
\end{equation}
This, together with (\ref{chain}), gives us the following integration by parts formula for $\omega\in\drh{M}{k},\mu\in\drh{M}{n-k-1}$:
\begin{equation}\label{intPart}
    \int_M d\omega\wedge\mu+(-1)^k\int_M\omega\wedge d\mu =\int_{\partial M} \omega\wedge\mu.
\end{equation}
Notice that, for now, there is no problem for the terms on the right in the last two formulas since it is always possible to restrict a smooth differential form to the boundary by pulling it back via the inclusion map $\partial M\hookrightarrow M$.\\
We define the space of $L^2$ differential forms of degree $k$ as the closure of $\drh{M}{k}$ w.r.t. the norm induced by (\ref{l2p}).\\
It is straightforward to check that this space is given by differential forms of degree $k$ whose coefficients are $L^2$ functions on any coordinate patch. Moreover, the wedge product (\ref{wedge}) extends naturally to forms with $L^2$ coefficients and retains property (\ref{swap}).\\
On $L^2\drh{M}{k}$, we may consider the external differential $d$ and codifferential $\delta$ (a priori only defined on $\drh{M}{k}$) as unbounded operators
\begin{equation}\label{ExtDiff}
    d\colon \mathcal{D}(d)\subset L^2\drh{M}{k}\to L^2\drh{M}{k+1}
\end{equation}
\begin{equation}\label{CoDiff}
    \delta\colon \mathcal{D}(\delta)\subset L^2\drh{M}{k}\to L^2\drh{M}{k-1}.
\end{equation}
In order to find fitting domains for these maps, we introduce the following spaces.

\begin{defn}[Weak external differential]
    For a differential form $\omega\in L^2\drh{M}{k}$ we say that $\mu\in L^2\drh{M}{k+1}$ is the weak external differential of $\omega$ if, for every $\sigma\in C^{\infty}_0\drh{M}{k+1}$,
        \begin{equation}\label{wexd}
        \scprod{\omega}{\delta \sigma}_{L^2\drh{M}{k}}=\scprod{\mu}{\sigma}_{L^2\drh{M}{k+1}}.
        \end{equation}
    where $C^{\infty}_0\drh{M}{k+1}$ is the space of smooth differential forms vanishing on the boundary.\\
    We denote the space of forms with weak external differential in $L^2\drh{M}{k+1}$ as $H\drh{M}{k}$ and we regard it as the domain of the unbounded operator (\ref{ExtDiff}).\\
\end{defn}
\begin{defn}[Weak codifferential]\label{codi}
    For a differential form $\mu\in L^2\drh{M}{k}$ we say that $\omega\in L^2\drh{M}{k-1}$ is the weak codifferential of $\mu$ if, for every $\sigma\in C^{\infty}_0\drh{M}{k-1}$,
        \begin{equation}\label{wcod}
        \scprod{\mu}{d\sigma}_{L^2\drh{M}{k}}=\scprod{\omega}{\sigma}_{L^2\drh{M}{k-1}}.
        \end{equation}
    We denote the space of forms with weak codifferential in $L^2\drh{M}{k-1}$ as $H^{\star}\drh{M}{k}$ and we regard it as the domain of the unbounded operator (\ref{CoDiff}).
\end{defn}
\begin{prop}\label{closeddom}
    The operators $d$ and $\delta$, on the respective domains $H\drh{M}{k}$ and $H^{\star}\drh{M}{k}$, are closed.
\begin{pf}
    Consider a sequence $\{\omega_n\}$ in $H\drh{M}{k}$.\\
    Suppose $\omega_n\rightarrow\omega$ strongly on $L^2\drh{M}{k}$ and $d\omega_n\rightarrow \tau$ strongly in $L^2\drh{M}{k+1}$.\\
    For every $\beta\in C^{\infty}_0\drh{M}{k+1}$, by the definition of weak external differential, we can write
\begin{align*}
  \scprod{\omega}{\delta\beta}=\lim_{n\rightarrow \infty}\scprod{\omega_n}{\delta \beta} = \lim_{n\rightarrow\infty} \scprod{d\omega_n}{\beta}=\scprod{\tau}{\beta}.
\end{align*}
Hence $d\omega=\tau$ weakly, and $\omega\in H\drh{M}{k}$.\\
A similar argument works for $\delta$.\hfill$\qed$
\end{pf}
\end{prop}
The spaces $H\drh{M}{k}$ and $H^{\star}\drh{M}{k}$ inherit the $L^2$ norm of $L^2\drh{M}{k}$, but also come naturally equipped with graph norms that make them Hilbert spaces, because of the previous proposition:
\begin{align*}
    \norm{\omega}_{H\drh{M}{k}}^2&\coloneqq \norm{\omega}_{L^2\drh{M}{k}}^2+\norm{d\omega}_{L^2\drh{M}{k+1}}^2,\\
    \norm{\mu}_{H^{\star}\drh{M}{k}}^2&\coloneqq \norm{\mu}_{L^2\drh{M}{k}}^2+\norm{\delta\mu}_{L^2\drh{M}{k-1}}^2.
\end{align*}
We may also regard the Hodge star operator as a map between $L^2$ forms:
\[
    \star\colon L^2\drh{M}{k}\rightarrow L^2\drh{M}{n-k}
\]
This is natural as the Hodge star is $C^{\infty}$-linear (that is $\star(f\omega)=f\star\omega$  for all $f\in C^{\infty}(M)$), hence extending coefficients to $L^2$ functions is always possible, and the property (\ref{invHS}) still holds.\\
By density of $\drh{M}{k}$ in $L^2\drh{M}{k}$, the Hodge star on $L^2$ forms is still an isometry with respect to the $L^2$ product.

\begin{prop}
    The weakly defined $d$ and $\delta$ coincide with their smooth counterparts on $\drh{M}{k}.$\\
    Also, (\ref{delta}) holds for the $L^2$ versions of $d,\delta$ and $\star$. Moreover, we have the following relations between the spaces of weakly differentiable forms:
    \begin{align*}
        \star \left(H\drh{M}{k}\right)&=H^{\star}\drh{M}{n-k}\\
        \star \left(H^{\star}\drh{M}{k}\right)&=H\drh{M}{n-k}
    \end{align*}
\begin{pf}
First, we prove that (\ref{wexd}) holds for $\omega\in \drh{M}{k}$ and $\mu=d\omega$ in the smooth sense. By using (\ref{intPart}) we indeed get
\begin{align*}
    \MoveEqLeft
    \scprod{d\omega}{\sigma}_{L^2\drh{M}{k+1}}=\int_Md\omega\wedge\star\sigma =(-1)^{k-1}\int_M\omega\wedge d\star\sigma\\
    &= (-1)^{k-1}\scprod{\omega}{\star^{-1}d\star\sigma}_{L^2\drh{M}{k}}=\scprod{\omega}{\delta\sigma}_{L^2\drh{M}{k}}
\end{align*}
for all $\sigma\in C^{\infty}_0\drh{M}{k+1}$.\\
An analogous calculation shows that (\ref{wcod}) holds for $\mu\in\drh{M}{k}$ and $\omega=\delta\mu$ is the smooth sense.\\
It follows that (\ref{delta}) still holds on weak forms, since the two sides of the equality coincide on $\drh{M}{k}$ which is dense in $H^{\star}\drh{M}{k}$.\\
Let now $\omega\in H\drh{M}{k}$, we may compute the codifferential of $\star\omega$ just by using the definition of $\delta$ (\ref{delta}) and (\ref{invHS}):
\[
    \delta\star\omega=(-1)^k\star^{-1}d\star\star\omega=(-1)^{k+(n-k)k}\star^{-1}d\omega.
\]
The right side of the last equation is an element of $L^2\drh{M}{n-k-1}$, hence $\star\omega\in H^{\star}\drh{M}{n-k}$\\
Similarly, for $\mu\in H^{\star}\drh{M}{k}$, we may compute the external differential of $\star\mu$ as
\[
    d\star\mu=\star\star^{-1}d\star\mu=(-1)^k\star\delta\mu.
\]
Again, the last term is an element of $L^2\drh{M}{n-k+1}$, thus $\star\delta\in H\drh{M}{n-k}$.\\
We have shown 
\begin{align*}
    &\star \left(H\drh{M}{k}\right)\subset H^{\star}\drh{M}{n-k}\\
    &\star \left(H^{\star}\drh{M}{k}\right)\subset H\drh{M}{n-k},
\end{align*}
and since, by (\ref{invHS}), applying the Hodge star twice only amounts to a sign change, we have proven our claim. \hfill $\qed$
\end{pf}
\end{prop}
It is possible to define the concept of distributions on a manifold due to their local nature, this leads to defining Sobolev spaces $H^s(M)$, for which we refer to Chapter 4 in \cite{PD1:96}. On a coordinate patch, this concept corresponds to the usual Sobolev spaces on an open domain in $\R^n_+$.\\
In particular we will be interested in $H^1(M)$ and $H^{1/2}(M)$.\\
We are also interested in differential forms whose coefficients, when picking a coordinate patch, are in these spaces, which we will denote, respectively, as
\[
    H^1\drh{M}{k},\quad H^{1/2}\drh{M}{k}.
\]
Note also that, since the coefficients of a form in $H^1\drh{M}{k}$ on a coordinate patch have well-defined weak partial derivatives, we have the inclusions
\[ H^1\drh{M}{k}\subset H\drh{M}{k}\quad H^1\drh{M}{k}\subset H^{\star}\drh{M}{k}.\]
These are, in almost all cases, strict inclusions, with the notable exceptions
\begin{align*}
    H\drh{M}{0}=H^1\drh{M}{0},\\
    H^{\star}\drh{M}{n}=H^1\drh{M}{n},
\end{align*}
which hold since, on a coordinate patch, for the space of $0$-forms (resp. $n$-forms), the operator $d$ (resp. $\delta$), in coordinates, is just $\nabla$ (resp. $(-1)^n\nabla$).
\subsection{Traces and an extension of Stokes theorem}
Just like for domains in $\R^n$, we may trace a $H^1$ form onto the boundary:
\[
    \tr\colon H^1\drh{M}{k}\to H^{1/2}\drh{\partial M}{k}
\]
This trace map is surjective and can be easily constructed by working on a coordinate patch and then using a partition of unity on $M$ to glue the local trace maps together (see Chapter $4$ in \cite{PD1:96}).\\
For a smooth form $\omega\in\drh{M}{k}$, $\tr\omega\in\drh{\partial M}{k}$ is just, once again, the pull-back on the boundary.\\
By density of $\drh{M}{k}$ in $H^1\drh{M}{k}$, this allows us to rewrite the formula (\ref{intPart}) for $\omega\in H^1\drh{M}{k-1}$ and $\mu\in H^1\drh{M}{n-k}$
\begin{equation}\label{stH1}
    \int_M d\omega\wedge \mu + (-1)^k\int_M \omega\wedge d\mu =\int_{\partial M} \tr \omega\wedge \tr \mu 
\end{equation}
The trace map is surjective, and it endows the space $H^{\frac{1}{2}}\drh{\partial M}{k}$ with the range norm
\begin{equation}\label{halfnorm}
    \norm{\tau}_{H^{1/2}\drh{\partial M}{k}}\coloneqq \inf_{\tr\omega=\tau}\norm{\omega}_{H\drh{M}{k}},
\end{equation}
which, trivially, makes $\tr$ continuous w.r.t the $H\drh{M}{k}$ norm.\\
In order for this norm to be well-defined, we need the folowing.
\begin{lem}\label{closedker}
    The kernel of $\tr$ is closed with respect to the $H\drh{M}{k}$ norm.
    \begin{pf}
        Consider a sequence $\{\omega_n\}$ in $\operatorname{ker}\tr$, converging to $\omega\in H^1\drh{M}{k}$ with respect to the $H\drh{M}{k}$ norm. For any $\theta\in H^1\drh{M}{n-k-1}$ we may write, using (\ref{stH1}),
        \begin{align*}
            \MoveEqLeft
            \abs*{\int_{\partial M} \tr(\omega-\omega_n)\wedge\tr\theta}\\
            &\leq \abs*{\int_M d(\omega-\omega_n)\wedge\theta} + \abs*{\int_M (\omega-\omega_n)\wedge d\theta}\\
            &= \abs{\scprod{d(\omega-\omega_n)}{\star^{-1}\theta}} +\abs{\scprod{\omega-\omega_n}{\star^{-1}d\theta}}\\
            &\leq \norm{d(\omega-\omega_n)}_{L^2\drh{M}{k+1}}\norm{\theta}_{L^2\drh{M}{n-k-1}} \\
            &\quad +\norm{\omega-\omega_n}_{L^2\drh{M}{k}}\norm{d\theta}_{L^2\drh{M}{n-k}}\\ 
            &\leq \norm{\omega-\omega_n}_{H\drh{M}{k}}\norm{\theta}_{H\drh{M}{n-k-1}},
        \end{align*}
        where for the second inequality we have used Chauchy-Schwarz and the fact that $\star$ is an isometry.\\
        Since $\omega-\omega_n$ converges to $0$ with respect to the $H\drh{M}{k}$ norm, we obtain that for all $\theta\in H^1\drh{M}{n-k-1}$
        \[
            \int_{\partial M} \tr\omega\wedge\tr\theta=\int_{\partial M} \tr(\omega-\omega_n)\wedge\tr\theta = 0.
        \]
        Since the space $H^{1/2}\drh{\partial M}{n-k-1}$ contains the smooth forms of degree $k$ on $\partial M$, and the trace map is surjective, this condition implies that $\tr\omega=0$.\hfill $\qed$
    \end{pf}
\end{lem}
We denote $H^1_0\drh{M}{k}\subset H^1\drh{M}{k}$ as the space of differential forms of degree $k$ with zero trace.\\
We can now use  (\ref{stH1}) to extend $\tr$ to $H^{\star}\drh{M}{k}$. Indeed, for $\omega\in H^1\drh{M}{k-1}$ and $\mu\in H^1\drh{M}{n-k}$, we have a similar calculation as in the previous proof:
\begin{align*}
    \MoveEqLeft
    \abs*{\int_{\partial M} \tr\omega\wedge\tr\mu}\\
    &\leq
    \abs*{\int_M d\omega\wedge \mu}+\abs*{\int_M \omega\wedge d\mu}\\
    &\leq \norm{d\omega}_{L^2\drh{M}{k}}\norm{\mu}_{L^2\drh{M}{n-k}} \\
    &\quad + \norm{\omega}_{L^2\drh{M}{k-1}}\norm{d\mu}_{L^2\drh{M}{n-k+1}}\\
    &\leq \norm{\omega}_{H\drh{M}{k-1}}\norm{\mu}_{H\drh{M}{n-k}}.
\end{align*}
Hence for $\tau\in H^{1/2}\drh{\partial M}{k-1}$, by taking the infimum over all $\omega\in H^1\drh{M}{k-1}$ such that $\tr\omega=\tau$, we get
\begin{equation}\label{con}
    \abs{\scprod{\tau}{\star^{-1}\tr\mu}_{L^2\drh{\partial M}{k-1}}}\leq \norm{\tau}_{H^{1/2}\drh{\partial M}{k-1}}\norm{\mu}_{H\drh{M}{n-k}}.
\end{equation}
Please note that the $\star^{-1}$ appearing on the left is the inverse of the Hodge star on the boundary.\\
This formula shows that the map
\[
    \scprod{\ \cdot\ }{\star^{-1}\tr\mu}:H^{1/2}\drh{\partial M}{k-1}\to \R
\]
is continuous for any choice of $\mu\in H^1\drh{M}{n-k}$, hence it is an element of the dual space of $H^{1/2}\drh{\partial M}{k-1}$, which we will denote as $H^{-1/2}\drh{\partial M}{k-1}$, where $L^2\drh{\partial M}{k-1}$ is the pivot space between the two.\\
Formula (\ref{con}) also tells us the linear map $\star^{-1}\tr$ extends by density to a continuous map
\[
    \star^{-1}\tr\colon H\drh{M}{n-k}\to H^{-1/2}\drh{\partial M}{k-1}.
\]
Finally, after pre-composing with the Hodge star on $M$, we have the desired extension of the trace
\begin{equation}\label{star}
    \tr^{\star}\coloneqq\star^{-1}\tr\star\colon H^{\star}\drh{M}{k}\to H^{-1/2}\drh{\partial M}{k-1}
\end{equation}
Using the notation we introduced, (\ref{stH1}) may now be rewritten more neatly as
\begin{equation}\label{pt}
    \scprod{d\omega}{\mu}_{L^2\drh{M}{k+1}}-\scprod{\omega}{\delta\mu}_{L^2\drh{M}{k}}=\scprod{\tr\omega}{\tr^{\star}\mu}_{H^{1/2},H^{-1/2}}
\end{equation}
for $\omega\in H^1\drh{M}{k}$ and $\mu\in H^{\star}\drh{M}{k+1}$, where on the right we have the duality pairing between $H^{1/2}\drh{\partial M}{k}$ and $H^{-1/2}\drh{\partial M}{k}$, with $L^2\drh{\partial M}{k}$ being the pivot space.\\
Formula (\ref{pt}) is obtained from (\ref{stH1}) by a straightforward computation (first assuming $\omega$ and $\mu$ to be $H^1$ and then extending by density)
\begin{align*}
    \MoveEqLeft
    \scprod{d\omega}{\mu}_{L^2\drh{M}{k+1}}-\scprod{\omega}{\delta\mu}_{L^2\drh{M}{k}}\\
    &=\int_M d\omega\wedge\star\mu -\int_M \omega\wedge \star\delta\mu \\
    &=\int_M d\omega\wedge \star\mu +(-1)^k\int_M\omega\wedge d\star\mu  \\
    &=\int_{\partial M} \tr\omega\wedge\tr\star\mu \\
    &=\scprod{\tr\omega}{\tr^{\star}\mu}_{H^{1/2},H^{-1/2}}.
\end{align*}
We end this section by proving that the new trace we defined is surjective
\begin{prop}\label{surj}
    The map (\ref{star}) is onto.
    \begin{pf}
        Consider a $\tau\in H^{-1/2}\drh{\partial M}{k-1}$, for every $\omega\in H^1\drh{M}{k-1}$ the map
        \begin{equation}\label{this}
            F_1(\omega)\coloneqq\omega \mapsto \scprod{\tr\omega}{\tau}_{H^{1/2},H^{-1/2}}
        \end{equation}
        is well-defined, and it is continuous on $H^1\drh{M}{k-1}$ with respect to the $H\drh{M}{k-1}$ norm, since the trace is continuous (because of our choice of the norm on $H^{1/2}\drh{\partial M}{k-1}$ (\ref{halfnorm})).\\
        Since $H^1\drh{M}{k-1}$ is dense in $H\drh{M}{k-1}$, $F_1$ can be extended to a linear functional $F$ on $H\drh{M}{k-1}$
        \[
            F\colon H\drh{M}{k-1}\to \R.
        \]
        Since $H\drh{M}{k-1}$ is a Hilbert space, by the Riesz representation theorem there exists $\beta\in H\drh{M}{k-1}$ such that
        \begin{align}\label{wow}
            \MoveEqLeft
            \scprod{\omega}{\beta}_{L^2\drh{M}{k-1}}+\scprod{d\omega}{d\beta}_{L^2\drh{M}{k}}=\nonumber \\
            &=\scprod{\omega}{\beta}_{H\drh{M}{k-1}}=F(\omega)
        \end{align}
        for every $\omega\in H\drh{M}{k-1}$.\\
        In particular for every $\omega\in C^{\infty}_0\drh{M}{k-1}\subset H^1\drh{M}{k-1}$ we obtain
        \begin{align*}
            \MoveEqLeft
            \scprod{\omega}{\beta}_{L^2\drh{M}{k-1}}+\scprod{d\omega}{d\beta}_{L^2\drh{M}{k}}=F(\omega)\\
            &= F_1(\omega)=\scprod{\tr\omega}{\tau}_{H^{1/2},H^{-1/2}}=0.
        \end{align*}
        That is, for all $\omega\in C^{\infty}_0\drh{M}{k-1}$ we have
        \[
            \scprod{d\beta}{d\omega}_{L^2\drh{M}{k}}=-\scprod{\beta}{\omega}_{L^2\drh{M}{k-1}}.
        \]
        Confronting this with Definition \ref{codi}, we have that
        \[
            \delta d\beta=-\beta
        \]
        in the weak sense. Thus, defining $\mu\coloneqq d\beta\in H^{\star}\drh{M}{k}$, we can rewrite (\ref{wow}) as 
        \[
            \scprod{d\omega}{\mu}_{L^2\drh{M}{k}}-\scprod{\omega}{\delta\mu}_{L^2\drh{M}{k-1}}= \scprod{\tr\omega}{\tau}_{H^{1/2},H^{-1/2}}
        \]
        for all $\omega\in H^1\drh{M}{k-1}$.\\
        Combining this with (\ref{pt}), we have that 
        \[
            \scprod{\tr\omega}{\tr^{\star}\mu}_{H^{1/2},H^{-1/2}}=\scprod{\tr\omega}{\tau}_{H^{1/2},H^{-1/2}}
        \]
        for every $\omega\in H^1\drh{M}{k-1}$.\\
        Since the trace is surjective, $\tr^{\star}\mu$ and $\tau$ coincide on the whole $H^{1/2}\drh{\partial M}{k-1}$, hence $\tr^{\star}\mu=\tau$.
        \hfill $\qed$
    \end{pf}
\end{prop}

\section{Boundary spaces and well-posedness}
In this section, we briefly recall the classical approach used in port-Hamiltonian theory for proving the well-posedness of a system. This approach is fundamentally rooted in the functional analytic theory as presented in Chapter 3 of \cite{Gor:91}, and it was first shown in the context of port-Hamiltonian systems in \cite{YZB:06} and further explored in \cite{KZ:15} for the case of multidimensional spacial variables.\\
In particular, we refer to Section 2 of \cite{KZ:15} for a comprehensive and modern exposition of the subject matter, which we will only briefly review.\\

For the Hilbert space $\mathcal{X}$ and the unbounded operator $A:\mathcal{D}(A)\subset \mathcal{X}\to \mathcal{X}$ we consider the abstract system
\begin{equation}\label{sys}
    \dot{x}(t)=Ax(t)\quad  \forall t>0,
\end{equation}
where $x(t)$ is a time-dependent function in $\mathcal{X}$.\\
We know that ($\ref{sys}$) is well-posed if $A$ generates a contraction semigroup. For a closed and densely defined operator $A$ this is equivalent to $A$ being maximally dissipative, by the Lumer-Phillips theorem (\cite{LPh:61}).
\begin{defn}[Maximally dissipative operator]
    A closed\\ and densely defined unbounded operator $A$ is called maximally dissipative if
    \begin{itemize}
        \item $\forall x\in \mathcal{D}(A)$ we have
        \[
            \scprod{Ax}{x}+\scprod{x}{Ax}=2\text{Re}\scprod{Ax}{x}\leq 0.
        \]
        \item For any other operator $B$ which has the first property, it holds
        \[ 
            B|_{\mathcal{D}(A)}=A\Rightarrow A=B.
        \]
    \end{itemize}
\end{defn}
The context in which we will apply the Lumer-Phillips result is the following. Suppose we have a skew-symmetric operator $A_0$, defined as a restriction of the operator $A$ in (\ref{sys}). We call $A_0$ the minimal operator and we aim to classify the possible extensions of $A_0$ which are maximally dissipative, and we do so by using the following fundamental tool.
\begin{defn}[Boundary triplet]\label{BouTr}
    Let $A_0$ be a skew-sym-\\ metric, unbounded, closed and densely defined linear operator on the Hilbert space $\mathcal{X}$. Then a boundary triplet for $A_0^*$ is given by 
    \begin{itemize}
        \item Another Hilbert space $\mathcal{B}$,
        \item Two bounded linear operators
        \[B_1:\mathcal{D}(A_0^*)\rightarrow \mathcal{B}\quad B_2:\mathcal{D}(A_0^*)\rightarrow \mathcal{B}'\]
        which we will sometimes refer to as the boundary operators,
    \end{itemize}
    such that the following properties hold:
    \begin{itemize}
        \item[i)] Abstract Green identity: For any $x,y\in \mathcal{D}(A_0^*)$ we have the formula
        \begin{equation}\label{AbsG}\scprod{A_0^*x}{y}+\scprod{x}{A_0^*y}=\scprod{B_1x}{B_2y}_{\mathcal{B},\mathcal{B}'}+\scprod{B_1y}{B_2x}_{\mathcal{B},\mathcal{B}'}\end{equation}
        \item[ii)] the component-wise map to the product space
        \[\begin{bmatrix}
            B_1 \\
            B_2
        \end{bmatrix}:\mathcal{D}(A_0^*)\longrightarrow \mathcal{B}\times\mathcal{B}'\]
        is onto.
    \end{itemize}
\end{defn}
We remark that different definitions of the concept of  a boundary triplet appear in literature, depending on which pairing one wants to use on the right side of (\ref{AbsG}); in our case the duality pairing between $\mathcal{B}$ and $\mathcal{B}'$ will suffice. In particular, we are interested in the case where this pairing is part of a Gelfand triple, with pivot space $\mathcal{B}_0$, that is
\begin{equation}\label{gelf}
    \mathcal{B}\subset\mathcal{B}_0\subset\mathcal{B}'.
\end{equation}
The following is Theorem $2.6$ from \cite{KZ:15} and it introduces a class of domains on which a restriction of the operator $A$ generates a contraction semigroup.
\begin{thm}\label{generation}
    Let $\mathcal{B}$ be a Hilbert space densely and continuously contained in a Hilbert space $\mathcal{B}_0$, let $\mathcal{B}'$ be the dual of $\mathcal{B}$ with pivot space $\mathcal{B}_0$, see (\ref{gelf}).\\
    Assume that $(\mathcal{B},B_1,B_2)$ is a boundary triplet for the operator $A_0^*$ on the Hilbert space $\mathcal{X}$.\\
    Let $[V_1\ V_2]\in \mathcal{L}(\mathcal{B}^2_0;\mathcal{K})$, where $\mathcal{K}$ is some Hilbert space, and define
    \[
        \mathcal{D}\coloneqq\left\{x\in\mathcal{D}(A_0^*)\ |\ B_2x\in\mathcal{B}_0\text{ and } \begin{bmatrix}
            V_1 & V_2
        \end{bmatrix}\begin{bmatrix}
            B_1\\
            B_2
        \end{bmatrix}x=0 \right\}.
    \]
    The following conditions are sufficient for the closure of $A_0^*|_{\mathcal{D}}$ to generate a contraction semigroup on $\mathcal{X}$:
    \begin{itemize}
        \item $\Real \scprod{a}{b}_{\mathcal{B}_0}\leq 0$ for all $a,b\in\mathcal{B}_0$ such that
        \[
            V_1a+V_2b=0,
        \]
        \item the operator inequality holds in $\mathcal{K}$:
        \[
            V_1V_2^*+V_2V_1^*\geq 0.
        \]
    \end{itemize}
\end{thm}

\section{The Hodge wave system}
\subsection{A boundary triplet for the Hodge wave equation}
For a compact, $n$-dimensional smooth Riemannian manifold $M$ with boundary, the Hodge Laplacian on $0$-forms (that is, on $L^2$ functions on the manifold) is defined as the unbounded operator
\[
    \Delta: L^2\drh{M}{0}\to L^2 \drh{M}{0},\quad \Delta=\delta d.
\]
The Hodge scalar wave equation on $M$ is then defined as
\begin{equation}\label{HW}
  u_{tt} + \Delta u = 0,
\end{equation}
where $u$ is a time-dependent differential form of degree $0$, that is $u:[0,T]\rightarrow L^2\drh{M}{0}= L^2(M)$.\\
Assuming $u(t)\in H\drh{M}{0}$ for $t\in [0,T]$, introducing the variables $\omega(t)=u_t(t)\in L^2\drh{M}{0}$ and $\nu(t)=du(t)\in L^2\drh{M}{1}$, we rewrite equation (\ref{HW}) as the system
\begin{equation}\label{HWsystem}
    \frac{\mathrm{d}}{\mathrm{d}t}\begin{bmatrix}
         \omega  \\
         \nu
    \end{bmatrix}=\begin{bmatrix}
        0 & -\delta \\
        d & 0  
    \end{bmatrix}
        \begin{bmatrix}
         \omega  \\
         \nu
    \end{bmatrix}.
\end{equation}
In order to study the well-posedness of this system, we consider $\mathcal{X}\coloneqq L^2\drh{M}{0}\times L^2\drh{M}{1}$ as our state space and the operator
\[
    L=  \begin{bmatrix}
        0 & -\delta \\
        d & 0  
        \end{bmatrix}
\]
as an unbounded operator on $\mathcal{X}$. 
Because of Proposition \ref{closeddom}, the maximal domain of $L$ is given by
\begin{equation}\label{dom}
    H\drh{M}{0}\times H^{\star}\drh{M}{1}.
\end{equation}
However this choice for the domain does not guarantee well-posedness for the problem as we have no boundary conditions for our equation.\\
In light of Section $3$, we seek a domain $\mathcal{D}$ such that $L|_{\mathcal{D}}$ is maximally dissipative, and thus generates a contraction semigroup on $\mathcal{X}$. Following the constructions in Section $3$, we consider the minimal operator
\[
    L_0\coloneqq -L|_{H_0\drh{M}{0}\times H_0\drh{M}{1}}.
\]
$L_0$ is skew-symmetric, in fact we have
\begin{align}\label{zbou}
\begin{aligned}
\MoveEqLeft
    \scprod{L_0(\omega_1,\mu_1)}{(\omega_2,\mu_2)}_{\mathcal{X}}+\scprod{(\omega_1,\mu_1)}{L_0(\omega_2,\mu_2)}_{\mathcal{X}} \\
    & = \scprod{\delta\mu_1}{\omega_2}_{L^2\drh{M}{0}}-\scprod{d\omega_1}{\mu_2}_{L^2\drh{M}{1}} \\
    &\quad + \scprod{\omega_1}{\delta\mu_2}_{L^2\drh{M}{0}}-\scprod{\mu_1}{d\omega_2}_{L^2\drh{M}{1}}=0
\end{aligned}
\end{align}
for every $(\omega_1,\mu_1),(\omega_2,\mu_2)\in H_0\drh{M}{0}\times H_0\drh{M}{1}$. Indeed the first and fourth term cancel each other out because of (\ref{pt}) and so do the second and the third.\\
Notice also that 
\[
    \mathcal{D}(L_0^*)=H\drh{M}{0}\times H^{\star}\drh{M}{1}.
\]
Indeed the calculation (\ref{zbou}) can be repeated for any $(\omega_2,\mu_2)\in H\drh{M}{0}\times H\drh{M}{1}$, substituting the second $L_0$ in the first line with the general operator $-L$.
\begin{thm}
    A boundary triplet $(\mathcal{B},B_1,B_2)$ for $L_0^*$ is given by
    \begin{align*}
    &\mathcal{B}=H^{1/2}\drh{\partial M}{0},\quad \mathcal{B}'=H^{-1/2}\drh{\partial M}{0},\\
    &B_1\colon \mathcal{D}(L_0^*)\to \mathcal{B},\quad \ B_1(\omega,\mu)=\tr\omega, \\ 
    &B_2\colon \mathcal{D}(L_0^*)\to\mathcal{B}',\quad B_2(\omega,\mu)=\tr^{\star}\mu.\end{align*}
    \begin{pf}
        We need to prove both properties in Definition \ref{BouTr}.\\
        A similar calculation to (\ref{zbou}), once again using (\ref{pt}), proves the abstract Green identity:
        \begin{align*}
        \MoveEqLeft
        \scprod{L_0^*(\omega_1,\mu_1)}{(\omega_2,\mu_2)}_{\mathcal{X}}+\scprod{(\omega_1,\mu_1)}{L_0^*(\omega_2,\mu_2)}_{\mathcal{X}}=\\
        & = \scprod{-\delta\mu_1}{\omega_2}_{L^2\drh{M}{0}}+\scprod{d\omega_1}{\mu_2}_{L^2\drh{M}{1}} \\
        &\quad + \scprod{\omega_1}{-\delta\mu_2}_{L^2\drh{M}{0}}+ \scprod{\mu_1}{d\omega_2}_{L^2\drh{M}{1}}\\
        &= \scprod{\tr\omega_2}{\tr^{\star}\mu_1}_{H^{1/2},H^{-1/2}}
        +\scprod{\tr\omega_1}{\tr^{\star}\mu_2}_{H^{1/2},H^{-1/2}}.
        \end{align*}
        The surjectivity of the component-wise map follows trivially from the surjectivity of the trace operators $\tr$ and $\tr^{\star}$ (Proposition \ref{surj}). \hfill $\qed$
    \end{pf}
\end{thm}

Having constructed a boundary triplet for our operator $L$, we can now use Theorem \ref{generation} to find boundary conditions which guarantee well-posedness of the problem.

\begin{thm}\label{final}
    Let $V_1,V_2\in\mathcal{L}(L^2\drh{\partial M}{0},\mathcal{K})$, for some Hilbert space $\mathcal{K}$, be such that
    \begin{itemize}
        \item $\Real\scprod{\theta}{\sigma}_{L^2\drh{\partial M}{0}}\leq 0$ for every $\theta,\sigma\in L^2\drh{\partial M}{0}$ satisfying
        \[
            V_1\theta+V_2\sigma=0,
        \]
        \item the operator inequality holds in $\mathcal{K}$:
        \[ 
            V_1V_2^*+V_2V_1^*\geq0.
        \]
    \end{itemize}
    Then the operator $L|_{\mathcal{D}}$, where
    \begin{align*}
        \MoveEqLeft
        \mathcal{D}\coloneqq \{ (\omega,\mu)\in H\drh{M}{0}\times H^{\star}\drh{M}{1}\ |\\
        &\quad \tr^{\star}\mu\in L^2\drh{\partial M}{0},\quad V_1\tr\omega+V_2\tr^{\star}\mu =0\},
    \end{align*}
    generates a contraction semigroup on $\mathcal{X}$.
\end{thm}
We show two examples of applications of this theorem:
\begin{itemize}
    \item Choosing $V_1=\operatorname{id}$, $V_2=0$ and $\mathcal{K}=L^2\drh{\partial M}{0}$, the second condition in Theorem \ref{final} is trivially satisfied and the first one reduces to (remembering our original choice $\omega=u_t$)
    \[
        \tr u_t=\tr\omega=0.
    \]
    Thus this choice gives us well-posedness for Neumann boundary conditions on $u$.
    \item Choosing $V_1=0$, $V_2=\operatorname{id}$ and $\mathcal{K}=L^2\drh{\partial M}{0}$, we once again get the second condition for free and obtain
    \[
        \tr^{\star}(du)=\tr^{\star}\mu=0.
    \]
    It can be shown that this is equivalent to asking that the normal component of the gradient of $u$ on the boundary vanishes.
\end{itemize}
\subsection{Modelling an oscillating body}
    In applications where the manifold $M$ represents an oscillating object, we may want to introduce physical parameters like the mass density $\rho$ and the Young modulus $T$, which we will assume to be strictly positive $C^{\infty}$ functions. In this case, the Hodge wave equation can be rewritten as
    \begin{equation}\label{HWc}
        u_{tt}=-\rho^{-1}\delta(T du),
    \end{equation}
    and, after substituting $\omega(t)=\rho u_t(t)\in L^2\drh{M}{0}$ and $\nu(t)=du(t)\in L^2\drh{M}{1}$, the system (\ref{HWsystem}) becomes
    \begin{equation}\label{HWsysH}
    \frac{\mathrm{d}}{\mathrm{d}t}\begin{bmatrix}
         \omega  \\
         \nu
    \end{bmatrix}=\begin{bmatrix}
        0 & -\delta \\
        d & 0  
    \end{bmatrix}\begin{bmatrix}
        \rho^{-1} & 0\\
        0 & T
    \end{bmatrix}
        \begin{bmatrix}
         \omega  \\
         \nu
    \end{bmatrix}.
\end{equation}
    The matrix 
    \[
        \mathcal{H}\coloneqq
            \begin{bmatrix}
            \rho^{-1} & 0\\
            0 & T
            \end{bmatrix}
    \]
    is called the Hamiltonian density.\\
    Since $\mathcal{H}$ is coercive, the question of well-posedness of (\ref{HWsysH}) on the Hilbert space $\mathcal{X}_{\mathcal{H}}$, which is nothing but our original space $\mathcal{X}$ with a modified inner product:
    \[
        \scprod{x}{y}_{\mathcal{X}_\mathcal{H}}\coloneqq \scprod{x}{\mathcal{H}y}_{\mathcal{X}},
    \]
    is equivalent to that of the well-posedness of (\ref{HWsystem}) on $\mathcal{X}$, as proven in Theorem $7.2.3$ in \cite{JZ:12}.
\begin{ack}
    The author recognizes and thanks his PhD supervisor Prof Dr H. J. Zwart from the University of Twente (Enschede, NL) for the support during the writing and revision of the paper.
\end{ack}

\bibliography{ifacconf}


\end{document}